%% file: v1.tex
\documentclass[10pt,letterpaper]{IEEEtran}
\usepackage{amsmath,amssymb,amsfonts,amsthm}
\usepackage{mathbbol}
\usepackage{gensymb}
\usepackage{balance}
\usepackage{graphicx}
\usepackage{subfigure}

\usepackage{cite}
\usepackage[dvipsnames,usenames]{color}
\usepackage[colorlinks, linkcolor=Blue, filecolor =Red,citecolor=RoyalPurple]{hyperref}
\graphicspath{{./figures/}}
\usepackage{scalerel}
\usepackage{nopageno}
\IEEEoverridecommandlockouts
\usepackage[lmargin=0.76in,rmargin=0.76in,tmargin=0.73in,bmargin=0.82in]{geometry}

\def\powerRated{\ensuremath{p^{\text{rated}}}}
\def\qDisturb{\ensuremath{q^\text{d}}}

\def\COP{\ensuremath{\eta_{\text{COP}}}}

\def\COPchw{\ensuremath{\eta_{\text{COP}}^\text{ch}}}

\def\degK{\ensuremath{\degree\text{K}}}
\def\degF{\ensuremath{\degree\text{F}}}
\def\degC{\ensuremath{\degree\text{C}}}
\def\ratedpower{\ensuremath{p^{\mathrm{rated}}}}

\def\stpt{\ensuremath{\text{sp}}}
\def\peq{\ensuremath{p^\text{eq}}}

\def\rgrid{\ensuremath{r^{\text{BA}}}}

\def\R{\mathbb{R}}

\def\mca{\ensuremath{{\dot{m}^\mathrm{CA}}}}
\def\msa{\ensuremath{{\dot{m}^\mathrm{SA}}}}

\def\Tma{\ensuremath{{T^\mathrm{MA}}}}
\def\Tca{\ensuremath{{T^\mathrm{CA}}}}

\def\Wma{\ensuremath{W^\mathrm{MA}}}
\def\Wca{\ensuremath{W^\mathrm{CA}}}

\def\Wra{\ensuremath{W^\mathrm{RA}}}

\def\qcc{\ensuremath{q^\mathrm{cd}}}

\def\prh{\ensuremath{p^\mathrm{rh}}}
\def\pcc{\ensuremath{p^\mathrm{cd}}}

\def\pfan{\ensuremath{p^\mathrm{fan}}}

\def\rd#1{{\color{red}{#1}}}

\newlength{\noteWidth}
\long\def\notes#1{\ifinner
          {\footnotesize #1}
          \else
          \marginpar{\parbox[t]{\noteWidth}{\raggedright\footnotesize #1}}
      \fi\typeout{#1}}

\def\pb#1{\notes{pb: \rd{#1}      }}

\usepackage{xparse}
\ExplSyntaxOn

\ExplSyntaxOff

\ExplSyntaxOn
\cs_new:Nn \pb_prop_gset_bykeys:Nn
 {
  \prop_clear:N \l__pb_temp_prop
  \keys_set:nn { pb/propbykey } { #2 }
  \prop_gset_eq:NN #1 \l__pb_temp_prop
 }
\cs_generate_variant:Nn \pb_prop_gset_bykeys:Nn { c }

\keys_define:nn { pb/propbykey }
 {
  unknown .code:n = \prop_put:Nxn \l__pb_temp_prop { \l_keys_key_tl } { #1 }
 }

\cs_generate_variant:Nn \prop_put:Nnn { Nx }

\NewDocumentCommand{\setupcollaborator}{mm}
 {
  \prop_new:c { g_collaborator_#1_prop }
  \pb_prop_gset_bykeys:cn { g_collaborator_#1_prop } { #2 }
 }

\NewDocumentCommand{\selectcollaborator}{m}
 {
  \prop_map_inline:cn { g_collaborator_#1_prop }
   {
    \tl_set:cn { ##1 } { ##2 }
   }
 }
\ExplSyntaxOff

\setupcollaborator{PBmacbookPersonal}
 {
  NAME=Prabir Mac Laptop,
  PBbibPATH=/Users/prabirbarooah/bibfilesFromDropbox/DiCElab_bib,
 }
\setupcollaborator{ACwindows}
 {
  NAME=AC at lab,
  PBbibPATH=../../../bibFiles/dicelab_bib,
  ACbibPATH=../../../bibFiles/austin_bib,
  JBbibPATH=../../../bibFiles/bibs,
  REFRESHMENT=ice cream,
 }

 \selectcollaborator{PBmacbookPersonal}

\begin{document}	

\title{Comments on characterizing demand flexibility to provide power grid services}

\author{
	\IEEEauthorblockN{Prabir Barooah\IEEEauthorrefmark{1}} 
	\thanks{\IEEEauthorrefmark{1} email:
          pbarooah@iitg.ac.in. Dept. of Electronics and Electrical
          Engineering, Indian Institute of Technology, Guwahati,
           India.}
}



\maketitle
\thispagestyle{empty}
\begin{abstract}
Many loads have flexibility in demand that can be used to provide ancillary services to power grids. A large body of literature exists on designing algorithms to coordinate actions of many loads to provide such a service. The topic of characterizing the flexibility of one or a collection of loads - to determine what kinds of demand deviation from the baseline is feasible - has also been studied. However, there is a large diversity in definitions of flexibility and methods proposed to characterize flexibility. As a result there are several gaps in the literature on flexibility characterization. Some approaches on flexibility characterization are based on ad-hoc approximations that lead to highly conservative estimates. In this paper we point out some of these issues and their implications, with the hope to encourage additional research to address them. 
\end{abstract}

\section{Introduction}\label{sec:intro}
It is now widely believed that inherent flexibility of demand that many loads enjoy can be used to provide useful ancillary services to power grids. Intentional change in demand to help the grid is referred to by various names, such as \emph{demand dispatch}~\cite{BrooksDemandPEM:10}, which can be used to provide a multitude of services, from peak demand reduction to frequency regulation to renewable generation following. Maintaining consumers' quality of service (QoS) while providing ancillary services is essential. Monthly energy bill is one such QoS, which can be maintained by keeping the long term energy consumption (kWh) of the loads unchanged. Rather, their power demand (kW) can be changed up and down from the nominal value (the so-called \emph{baseline} demand) so that consumers' QoS is not adversely affected. This type of demand dispatch is called \emph{virtual energy storage} (VES), since they provide the same service as a battery~\cite{bar:2019}.  

The literature shows extensive work on designing coordination algorithms so that a collection of loads can provide some ancillary service. Less effort has been spent on quantifying the flexibility, whether for one load or for a collection of loads. A precise characterization of demand flexibility of loads is just as important as algorithms to utilize flexibility. Otherwise loads might be asked to do too much and forced to violate their quality of service. Consumers may not adopt such technology. Or, the grid operators - the balancing authorities (BAs)- may not be able to plan appropriately and utilize demand flexibility effectively.

There is a diversity in definitions of demand flexibility, techniques to characterize/compute flexibility, intended grid service, and consumers' QoS constraints. Also, some topics have attracted significant attention, such as flexibility characterization of a collection of thermostatically controlled loads, which are typically small residential loads. But large commercial or industrial loads have not attracted as much attention. Similarly, questions on characterizing capacity of virtual energy storage using terminology used for real batteries (MW/MWh) have not been adequately debated.


In this paper we discuss some of these gaps regarding demand flexibility characterization. Many questions arise in discussing demand flexibility~\cite{LiTenBE:2022}. We do not attempt to address all. Our focus is on demand flexibility definitions and characterization methods, especially for heating, ventilation and air conditioning (HVAC) systems. Reviewing all the available work is not possible here. Rather, we try to identify the remaining gaps that need to filled to make the path for technology adoption smoother.




The rest of the paper is organized as follows. Section~\ref{sec:QoS} summarizes the consumers' constraints that limit demand flexibility from various loads. Section~\ref{sec:one} discusses demand flexibility definitions for a single load, while Section~\ref{sec:many} discusses flexibility of a collection of loads, both dictated by consumers' QoS. Section~\ref{sec:industry} discusses issues that are dictated by grid balancing authorities QoS constraints. Section~\ref{sec:conclusion} concludes the paper.

\section{QoS constraints of consumers and grid balancing authorities}\label{sec:QoS}
Consumers expect certain quality of service (QoS) from their appliances (loads). Grid operators or balancing authorities (BAs) also expect certain QoS from any technology that provides ancillary services. Consumers' QoS expectations vary
depending on the type of the load. For air conditioners it is indoor
temperature and humidity, while for water heaters it is
availability of sufficient hot water when needed. Apart from (i) space temperature, (ii) humidity, and (iii) indoor air quality, there are
additional QoS constraints when it comes to HVAC loads: (iv)
equipment lifetime, (v) noise, and (vi) monthly utility bill. In fact, monthly energy bill is a constraint for any consumer load and will not be explicitly mentioned from now on. For an EV, a constraint is sufficient state of charge (SoC) before a trip.

We believe that VES technology is
likely to be accepted by consumers only if it ensures that there is \emph{no
noticeable change} in their QoS from their \emph{baseline} values while providing VES service. Baseline refers to the scenario
when equipment only serve the needs of the consumers and no consideration is made for the needs of the power grid. In other words, the amount of   demand flexibility available is limited by QoS constraints.

\ifx 0
For residential air conditioners, demand flexibility comes from \pb{this can be shortened}
the fact that a small variation of the indoor temperature and humidity is 
not noticed by occupants. In fact, even when a thermostat is set at a
specific setpoint, the on/off actuation of the compressor and the fan
does not maintain (and is not capable of maintaining) the indoor
temperature at that setpoint. Rather the temperature varies within a
low and high limit. Unless the temperature exceeds this range, which
can occur due to small capacity of the equipment compared to the
cooling load, or due to transients (when the AC is first turned on
while the indoor is stil quite warm), consumers are happy.

Figure~\ref{fig:3444data} shows measurements of
indoor temperature inside a house in Gainesville, Florida, USA with a
conventional thermostat set at $74$ \degF. The air temperature
fluctuates approximately within $[71.5, 74]$\degF, so the lower limit
of the thermostat logic is probably close to 71.5. The lower limit is
not known to the occupant; it is hardcoded into the thermostat's
software. Yet, almost all occupants are OK with this level
of variation. Thus, any control algorithm to provides DD must keep the indoor
temperature within such a range. The question of humidity will be
addressed in Section~\ref{}.

\begin{figure}[ht]
  \centering
  \includegraphics[width=0.7\columnwidth]{myHouseTemperature_inF.pdf}
  \includegraphics[width=0.28\columnwidth]{rig.pdf}
  \caption{(Left) Measured indoor temperature in the author's (former)
    residence with the air conditioner operated with a thermostat. The
  data was collected with an Arduino and a hobby grade temperature
  sensor located within 4 ft. of the thermostat (right) during Nov 8-9, 2020.}
  \label{fig:3444data}
\end{figure}
\fi

Let us consider a specific load and mathematically express its QoS constraints for use in the sequel. Let $q(t)$ be the vector of QoS signals at time $t$ and $Q(t)$ be the set in which $q(t)$ needs to lie to satisfy a consumer's QoS. For instance, if the only QoS requirement is that the temperature $\theta(t)$ and humidity ratio\footnote{The humidity ratio of a moist air volume is the ratio of the mass of water vapor to the mass of dry air. It is related to but is distinct from relative humidity, which is the ratio of partial pressure of water vapor in moist air to the saturation vapor pressure at the same temperature~\cite{ASHRAE_handbook_fund:17}.} $W(t)$   lie in a predetermined range $[\theta_{\min}, \theta_{\max}]$, $[W_{\min}, W_{\max}]$ then $q = [\theta, W] \in \R^2$ and $Q(t) = [\theta_{\min}, \theta_{\max}] \times [W_{\min}, W_{\max}]$. For air conditioners with on/off actuation, there is a lock out constraint to avoid compressor damage: if the command $u(t)$ changes from 0 to 1 or vice versa at time $t$, it cannot change again before $t+\tau_{lock}$, where $\tau_{lock}$ is the lock-out time. One can model the lock out constraint by introducing a new state variable, $s(t)$ that keeps a count of the number of changes of the control command in any direction during the last lock out interval: $s(t):=\int_{t-\tau_{lock}}^t |\dot{u}(v)|dv$ (here $|\dot{u}(v)|$ is a shifted Dirac-delta function since $u(v_0^-)=0$ and $\dot{u}(v_0) = 1$ if the device is turned on at $v_0$). Then, the lock out QoS requirement becomes $s(t)  \in [0, 1]$ for all $t$, which can be represented as $q(t) \in Q(t)$, with $q(t) =[\theta(t), W(t), s(t)]^T$. 

Just like the consumers, the balancing authority (BA) too has requirements on its quality of (ancillary) service provided by the loads. Suppose  $\rgrid(t)$ is the grid's desired power consumption from a giant battery, meaning if a battery (or a virtual battery) consumes $\rgrid(t)$ then the demand supply imbalance in the grid will remain 0. The balancing authority will of course use many resources including fast ramping generators, flywheels and batteries in addition to loads to supply the total $\rgrid$. The balancing authority  needs to know what part of $\rgrid(t)$ can be supplied by a VES system made of flexible loads with a predefined bound on tracking error. The BA also needs to know for long-term planning what fraction of its expected imbalance can be reliably met by flexible loads. This is where demand flexibility characterization becomes critical. 


\ifx 0 
\subsection{The dogma of price  elasticity \rd{remove?}}
\pb{also discuss the work by LBNL, Sila Kilicote, Vrettos, etc. : there is a conflict: why leave room for up/down flex from the optimal? Their answer was if we are paid enough to be inefficient. But isn't it inefficient from an economic point of view? A better answer is that (can I show it?) there are multiple local minimum if there is some leeway in temperature. So it is better to choose a baseline that has the same energy cost as any other optimal sol, but that has more leeway? If I cannot prove it, pick the optimal and find out the time intervals in which the temp is not at the lowest bound. In those times we can increase demand to reduce temp, and then decrease demand to bring temp back to stpt., but in a zero mean fashion. Then the total cost , energy , is the same (unless TOU pricing) but we created flexibility! Wait, this is a proof that there is no unique min. So I have to examine the opt problem to make sure it is not a strictly convex problem. Sum of power is a convex function, and if the constraint set is convex, done! Have to have a theorem showing more than one solution exists; but that follows from the previous construction. Question is: why would the temperature not hug one of the boundaries to reduce energy use?????}

A large number of works have used an alternate model of flexibility
that is based on an assumed price elasticity of consumers: if
cost of energy is changed consumers will change their demand and
tolerate the resulting loss of QoS; see e.g.~\cite{LiOptimalPESGM:2011,XX}. In other words, customers would
\emph{happily} tolerate some discomfort if they are sufficiently
compensated financially.

Price elasticity to not just electricity but various other
goods and services has been widely analyzed, but its effectiveness on behavior modification is unclear. A good example is smoking. After the damaging effects of smoking became clear, price of cigarettes was raised in the West with an aim to reduce smoking. Smoking did reduce, but an analysis of the data reveals many anomalies. ``...the real price of cigarettes in Britain was lower in 1990 than in 1965, but per-capita consumption was 20 per cent lower. What really changed, beginning in the 2000s, was the culture around smoking.''~\cite{BerginEconomicsAeon:2022}. Also, consumers that continue to smoke in the USA and Europe today come from poorer segments of society that have the least financial incentive to smoke. ``Today, in affluent neighborhoods in Britain smoking rates are under 10 per cent, whereas in some poor ones it’s 50 per cent.''~\cite{BerginEconomicsAeon:2022}.

Another example is the recent increase in the cost of fuel in the UK: if the dogma of price
elasticity were true, most people would silently reduce demand in
response to increase in energy cost; these increases would not have
created a political storm~\cite{}.  Therefore, we stay away from any notion of flexibility that relies on the unproven and risky assumption that consumers will \emph{happily tolerate noticeable changes} in their Quality of Service (QoS) for a financial reward. 
\fi

\section{Flexibility of one load (consumer focus)}\label{sec:one}
We consider a flexible load whose (real) electrical power demand is
denoted by $p(t)$. This demand is assumed \emph{continuously
  variable} within the range $[0, \powerRated]$. This is the case for a commercial HVAC system, in which power variation is actuated by varying the speed of a variable speed motor of a supply air fan~\cite{linbarmeymid:2015} or by changing chilled water temperature setpoint~\cite{SuDemonstration-I-STBE:2015}.

It follows from the discussion in Sec.~\ref{sec:QoS} that a natural definition of 
\emph{demand flexibility capacity of a flexible load is the set of possible demand trajectories that satisfy its QoS constraint}. For a time horizon $[t_0, t_f]$, this set is
\begin{align}
  \label{eq:flexCapacity-general}
\Omega_{t_0}^{t_f}:=\left\{ p(t)_{t_0}^{t_f} | q(t) \in  Q(t), \forall t \in
  [t_0,t_f],  \right\}, 
\end{align}
with $q(0) \in Q(0)$. The reason for limiting to a finite time horizon is purely technical;
choosing an infinite horizon will require us to specify signal spaces
that will be distracting at this stage. We will omit the time duration in the sequel, and refer to the set simply as $\Omega$. Note that this definition is not new; similar definition has been used implicitly in other works such as~\cite{adetola2018building}. Also, one can define the flexibility set~\eqref{eq:flexCapacity-general} in terms of the demand deviation $\tilde{p}:=p - p^b$, where $p^b$ is the baseline demand, instead of the demand $p$.

To vary the demand from the baseline, some control command $u$ must be manipulated that can affect the demand $p$. Doing so will potentially affect the QoS $q$ as well. To compute the set $\Omega$, one needs models $M_{p}: u\to p$ and $M_q: u \to q$.  These models will be necessarily dynamic since the underlying processes have memory. Note that the role of disturbance (e.g., weather for HVAC) is implicit in these models. 

Given a grid-level reference for the total VES requirement, $\rgrid$, one can project $\rgrid$ to $\Omega$ to determine a feasible demand trajectory, $r^*$, for the load that comes closest to what the grid needs. A control system can then be designed to make $p$ track $r$. Alternatively, one can simultaneously plan an optimal demand reference in $\Omega$ and the corresponding control command to realize that demand, by solving:
\begin{align}\label{eq:planning}
  \begin{split}
  (r^*,u^*) & = \arg \min_{p,u} \|\rgrid - p \|, \\
        &\text{ s.t.  } q \in Q, u \in U, q = M_q(u), p = M_p(u)   
  \end{split}
 \end{align}
The solution $r^*$ to~\eqref{eq:planning} provides the optimal reference that is within capacity of the load. Since solving~\eqref{eq:planning} requires forecast of $\rgrid$  for the time duration involved, which will have uncertainty, one can close the loop by using receding horizon control as updated forecasts become available.

\subsection{Why the definition matters} 
When the demand and QoS models $M_p,M_q$ are simple, computing feasible or optimal demand deviations within capacity by solving optimization problems is not difficult. That might be one of the reasons why more attention has been paid to collection of loads (which we will discuss in Sec.~\ref{sec:many}), since ensuring local QoS while the collection tracks a reference is challenging even if each load is simple. However, when models are complex, such as that for large commercial HVAC systems, computing the flexibility set $\Omega$ or planning a optimal demand within the flexibility set is quite challenging. That requires approximations. These approximations are sometimes made in an ad-hoc manner, leading to inaccurate capacity estimates. We give an example next. Surprisingly, the model used in this case is quite simple.

\subsubsection{A conservative estimation of flexibility}
Consider a simple HVAC system, in which indoor temperature is the only QoS signal. The flexibility capacity
set~\eqref{eq:flexCapacity-general} becomes
\begin{align}\label{eq:flexCapacity-HVAC-T-alone}
  \Omega & :=\left\{ p(t)_{t_0}^{t_f} | \theta(t) \in  [\theta_{\min},
  \theta_{\max}] \right\} 
\end{align}
It is common in the literature to use a simpler notion of flexibility
capacity in this scenario, as a lower and
upper and lower bound on the demand that maintains temperature within
its allowed range~\cite{OldewurtelTowardsCDC:2015,Lin_ACC2018,ChenAggregatingACC:2018,HuanAssessmentEnergy:2021,KhurramMethodologyPowerTech:2021}:
\begin{align}  \label{eq:flexCapacityBound}
  \begin{split}
  \hat{\Omega}: & = [\underline{p}(t), \bar{p}(t) ]_{t_0}^{t_f} \\
  & \text{ s.t} \;\;
\forall p \in [\underline{p}(t), \bar{p}(t) ]_{t_0}^{t_f}, \theta(t) \in
  [\theta_{\min}, \theta_{\max}]   
  \end{split}
\end{align}
Our claim is that \emph{the set $\hat{\Omega}$ is a small subset of the true
flexibility set $\Omega$. In other words, using $\hat{\Omega}$ leads to a gross
under-estimation of the true flexibility capacity of the load, and
thus should be avoided.}

We prove this claim by providing a specific example. Consider the following widely used and simple ``resistor-capacitor network'' model of an HVAC system providing cooling:
\begin{align}  \label{eq:RCnetwork}
  C\frac{d\theta}{dt} = -\frac{1}{R}(\theta(t) - \theta^a(t)) + \qDisturb(t) - \COP  p(t)
\end{align}
where $R$ (\degK/kW) and $C$ (kJ/\degK) are the thermal resistance of the structure and its thermal capacitance, $\theta^a(t)$ is the ambient (outdoor) dry bulb temperature, \COP\ is the coefficient of performance, and $\qDisturb$ is the  disturbances that includes heat gains from occupants and appliances, solar irradiance, etc. The total \emph{cooling load (kW-thermal)} is $-\frac{1}{R}(\theta(t) - \theta^a(t)) + \qDisturb(t)$ which must be removed by the HVAC system to maintain a steady temperature. 

\begin{figure}[ht]
  \centering
  \includegraphics[width=\columnwidth]{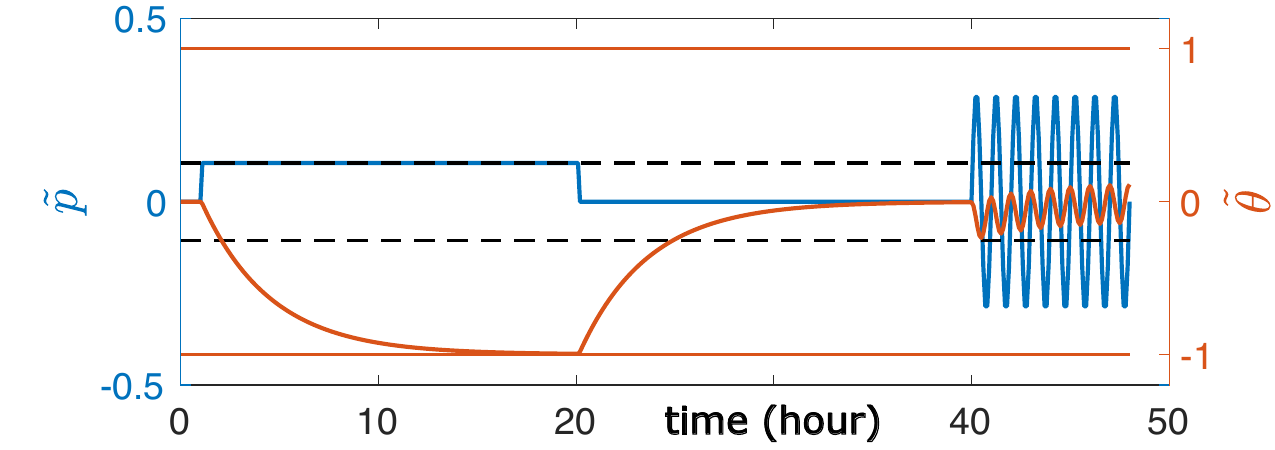}
  \caption{Simulation evidence: the largest envelope of power deviation is rather small ($\Delta_{\theta}/R\COP = 0.1055$ kW in this case; the dashed line) in order to ensure temperature constraints ($\pm 1$ \degC), while a sinusoidal variation with much larger amplitude, $0.3$ kW, can still maintain temperature within the same range if the frequency of the sinusoid is large enough (1 hour$^{-1}$ in this case).}
  \label{fig:inaccurateCapacity}
\end{figure}

The simplest and the most instructive case, for proving the claim, is  the time invariant one, when both the exogenous signals $\theta^a$  and $\qDisturb$ take constant values $\theta^a_0$ and $\qDisturb_0$, in which case an equilibrium
  the indoor temperature can be maintained at a constant setpoint
  $\theta^\stpt$ by demand $\peq$. A simple calculation from~\eqref{eq:RCnetwork} shows that   the equilibrium power demand is $\peq = \frac{1}{\COP}(\qDisturb_0 + \frac{1}{R}(\theta^a_0 - \theta^\stpt))$.  The ODE model of the temperature deviation $\tilde{\theta}(t):= \theta(t) - \theta^{\stpt}$ becomes
$    \dot{\tilde{\theta}}(t) = - \frac{1}{RC}\tilde{\theta}(t) - \frac{\COP}{C}\tilde{p}(t)$, where $\tilde{p}(t):=p(t) - \peq$. The transfer function from the electrical demand deviation to the
  temperature deviation is a low-pass filter: $  \displaystyle  G(s) =  \frac{- \frac{\COP}{C}}{s + \frac{1}{RC}}$. It follows from
  elementary frequency response of linear time invariant (LTI)
systems that if the allowed maximum temperature deviation is $\Delta_{\theta}$,
then the largest possible amplitude of a sinusoidal power deviation with frequency $\omega$ is $  A_{\max} = \frac{\Delta_{\theta}}{|G(j\omega)|} $. Due to the low pass nature of $G(j\omega)$, $A_{\max}$ is smallest at 0 frequency and its value increases as frequency increases. Thus, if the alternate definition of flexibility,~\eqref{eq:flexCapacityBound}, were used,
the largest envelope $\bar{p}-\underline{p}$ would be equal to $2\frac{\Delta_{\theta}}{|G(0)|} =  2\frac{\Delta_{\theta}}{R\COP}$. However this is not correct; a high frequency sinusoidal
component with a much higher amplitude still produces temperature
deviations that are small and does not violate QoS constraint of
the consumer, since the gain from power to temperature at that frequency is much smaller. Figure~\ref{fig:inaccurateCapacity} provides numerical verification. The parameters for the simulation are $\Delta_{\theta} = 1$ \degC, $\ratedpower - p^b_0 = \ratedpower - p^\text{eq} = 1$ kW, $R = 2.707$\degC/kW, $C = 1.283$ kWh/\degC, and $\COP=3.5$. The $R,C$ values are obtained by fitting the model to measurements from a real building and reported in~\cite{GuoAggregationEnB:2020}.  

Although the argument here is made in terms of a specific LTI model,
it clearly holds as long as the dynamic model from power deviation to
temperature deviation - whether linear or nonlinear - shows low-pass
characteristics. And, such low pass characteristics are expected due to the large thermal mass of buildings. Similarly, the argument above holds even with time-varying outdoor weather; the only
difference will be that the $\hat{\Omega}$ will be a time varying envelope. 

\ifx0
\rd{Three additional questions that we will discuss in the sequel are: (1)
Does humidity affect flexibility capacity? (2) What about flexibility
capacity of a collection of loads? (3) Can an HVAC system provide both demand flexibility and energy efficiency?}
\fi
\subsubsection{Definitions for deferrable loads}
Another line of flexibility definitions is based on the idea of
deferrable loads, such as pool pumps~\cite{meybarbusyueehr:TAC:2015} and - especially - electric vehicles
(EVs). These have a certain amount of energy demand over
a fixed time period but it can defer the power demand to some extent
as long as total energy demand is completed by the deadline.

A \emph{deferrable energy load} is defined
in~\cite{MadjidianEnergy_TPS:2018} as follows.
``a single deferrable energy load is characterized
by an arrival time $\tau \in \R$, an energy demand, $E$, a time
period, $T$, in which the demand must be filled, and a limit, $P$ on
its maximum power consumption. The energy consumed at time $t$ by a
load with arrival time $\tau$ is denoted by $e_{\tau}(t) =
\int_{-\infty}^{t}p_{\tau}(v)dv$, where $p_{\tau}$ is the
corresponding power consumption''. Any power consumption $0
\leq p(t) \leq \ratedpower$ is allowed as long as the energy
requirement is satisfies $e_{\tau}(t) = 0$ for $t \leq \tau$ and $
e_{\tau}(t) =  E$ for $t> \tau+T$.

An earlier set of papers also came up with similar
definitions but with different terminology; see~\cite{SubramanianRealtimeACC:2012,taxonomyPetersenACC:2013}. Ref.~\cite{taxonomyPetersenACC:2013} proposed a taxonomy
of demand flexibility as ``buckets, batteries, and bakeries''. A battery in that terminology is the same as
a deferrable energy load defined in~\cite{MadjidianEnergy_TPS:2018}, while a bakery is one in
which the energy must be consumed by a power trajectory
during a single interval. A bucket is similar to
a battery, but it does not have a predetermined energy demand that must be consumed by the deadline.

These definitions -  especially the deferrable energy load aka battery
- capture the flexibility of EVs, since most drivers want to charge
their EV batteries fully, which decides the parameter
$E$, by the beginning of every day, which decides $T$. 
But they are \emph{not} suitable for characterizing
an HVAC load's flexibility, although some works, such
as~\cite{SubramanianRealtimeACC:2012}, claim that they are. The energy requirement of an HVAC system
over any fixed period, say, 24 hours, is strongly determined by
outdoor weather. More importantly, constraints on QoS signals such as
indoor temperature are not part of the definition. It is rather
trivial to construct power demand signals that satisfy the
constraints in the definition a deferrable energy load but
fail to meet indoor temperature constraints. Let $T$ be 24 hours and $E$ be the baseline energy needed to maintain indoor temperature on a hot day. Running the air
conditioner on max power a cold day until $E$ is consumed, and then turning it off, will meet the requirements of the deferrable energy load definition but will violate indoor temperature constraint.

\emph{In contrast, the general defintion provided by~\eqref{eq:flexCapacity-general} can capture the constraints of deferrable loads/batteries with appropriate definitions of $q$ and $Q$.}

\subsection{(Under-appreciated) Role of humidity in HVAC flexibility}
HVAC systems are designed to maintain \emph{both temperature and humidity}, not just temperature.  Humidity is a key concern in many climate zones in
 the world, especially in the Southern and Western USA, and South East
 Asia. It is becoming a concern in areas that traditionally did not
 have to worry about humidity, such as parts of Europe, due to increase in extreme weather conditions. Yet, humidity is often ignored in the literature on demand flexibility. Depending on the type of equipment, this can lead to large errors in flexibility characterization. 

\begin{figure}[ht]
	\includegraphics[width=1.0\linewidth]{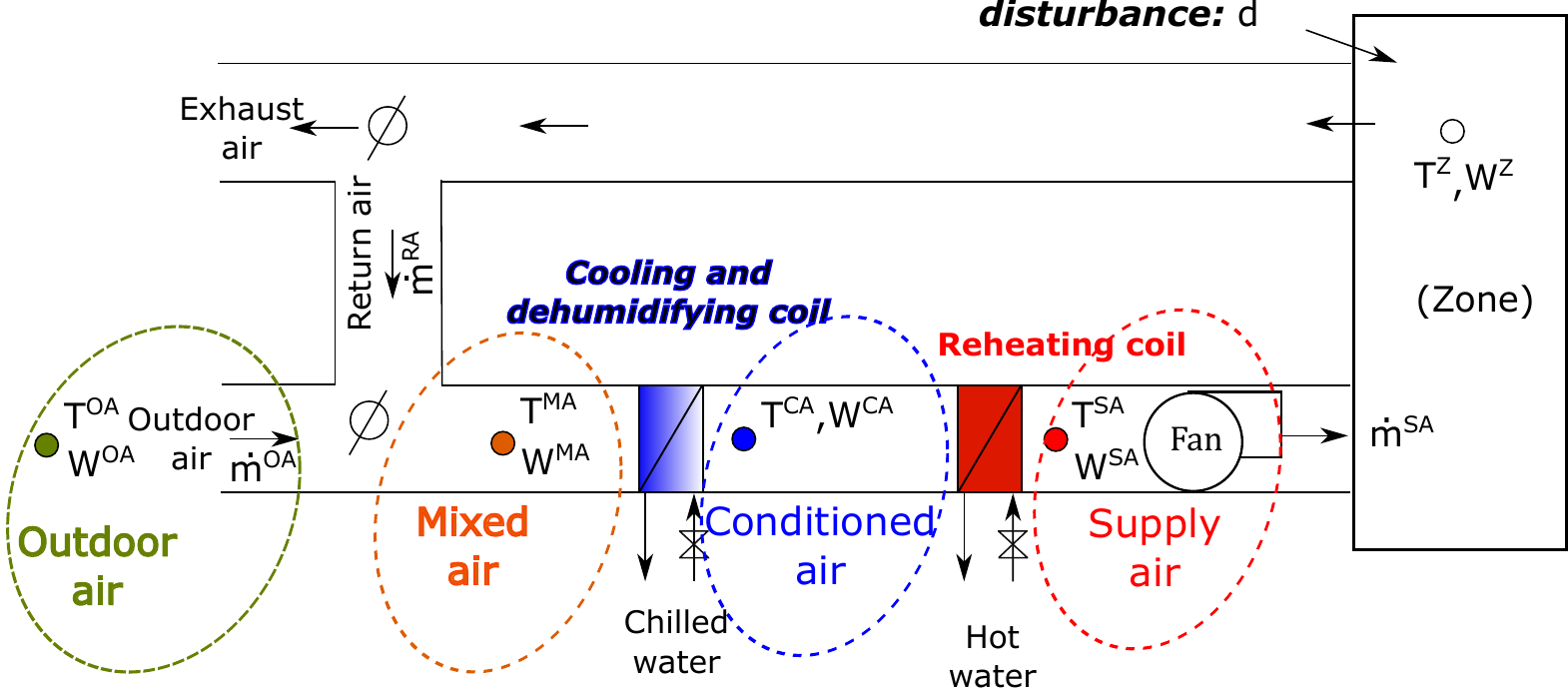}
	\caption{A single zone VAV HVAC system.}
	\label{fig:HVAC-single-zone}
      \end{figure}

To examine the effect of humidity on HVAC electric demand, let us examine a typical air handling unit in a commercial building, which is shown in
Figure~\ref{fig:HVAC-single-zone}. The subscripts MA, CA, RA, and OA in the figure refers to mixed air (before the cooling coil), conditioned air (after the cooling coil, right before delivery to the zones), return air (in same condition as  the building's interior), and outdoor air. The mixed air (MA) stream - mixture of outdoor air (OA) and return air (RA) - is cooled and dehumidified by passing over a cooling coil to produce conditioned air (CA) which is then supplied to the zones. The supply air (SA) to the zones is sometimes reheated since CA is often quite cold  (typically 55\degF).
Let us first look at the (thermal) power demand for cooling and dehumidification at an air handling unit (AHU), $\qcc$. This quantity can be approximated by the heat extracted from the mixed air stream by the chilled water in the cooling coil:
\begin{align}\label{eq:def-pcc}
\qcc & := \msa\left(h(\Tma,\Wma)-h(\Tca,\Wca)\right)
\end{align}
where $h(\cdot)$ is the \emph{specific enthalpy} of an air stream with (dry bulb) temperature $T$ and humidity ratio $W$:
\begin{align}
  \label{eq:sp-enth-def}
 h(T,W) := C_p T + W (h_g + C_{pw} T) 
\end{align}
where $C_p$ and $C_{pw}$ are the specific heat capacity of dry air and water ($1$ kJ/\degC/kg, $4.184$ kJ/kg/\degC, respectively),  and $h_g$ is the latent heat of evaporation of water (approx $2256$ kJ/kg at atmospheric pressure and 100 \degC)~\cite{ASHRAE_handbook_fund:17}. For temperatures encountered in HVAC systems, the second term within parenthesis is much smaller than the first, so we have $h(T,W) \approx C_p T + h_gW$. 

In the interest of simplicity, let us use the simplest possible model of electrical demand due to cooling and dehumidification, by dividing the thermal demand at the AHU, $\qcc$ with the effective COP of the chiller plant that produces the chilled water (discussed in Sec.~\ref{sec:largeHVAC}), yielding
\begin{align}
  \label{eq:p-coil-approx}
\pcc \approx  \frac{1}{\COPchw} (C_p (\Tma - \Tca) + h_g (\Wma - \Wca))  
\end{align}
Because of the need to maintain indoor humidity, the conditioned air -
downstream of the cooling coli - needs to have low moisture. This is
typically done by maintaining $\Tca$ near 55\degF, which is expected to maintain $\Wca \approx 0.004$ kg/kg, with the air downstream of the coil at $\approx 100\%$ relative humidity~\cite{WilliamsWhy:2013}. Since the ASHRAE
mandated target for comfortable indoor climate is 75\degF\ and 50\%
relative humidity, i.e. $\Wra = 0.009$ kg/kg, the return air can be expected to be at this condition. We first consider in detail the situation when
100\% of the return air is recirculated, so that the mixed air
conditions are the same as return air conditions. So, the specific
\emph{sensible cooling load} is
$C_p (\Tma - \Tca) \approx 20$ kJ/kg. Since $h_g = 2256$ kJ/kg, the \emph{latent cooling load}
$h_g(\Wma - \Wca) \approx 11$ kJ/kg. \emph{That is, the latent cooling load is similar to sensible cooling load in magnitude.} It follows from~\eqref{eq:p-coil-approx} that prediction of electricity demand for cooling with a  model that
ignores humidity can have a large error ($ \frac{11}{20+11} \approx 35\%$). In the more realistic situation - in which some outdoor air is brought in to maintain positive pressurization and indoor air quality - the situation can get worse in hot humid climate since the mixed air stream has outdoor air that is more humid than return air~\cite{fischer2003humidity}.

More importantly, \emph{the impact of humidity on HVAC electrical 
demand discussed above holds even when the climate or weather is not particularly humid}. The reason is that the return air has non-trivial amount of humidity  by design,
in interest of  occupant comfort, and since a large fraction of the mixed air is typically return air, the latent cooling load is still comparable to the sensible cooling load.  So the calculations done above are valid not only for
Gainesville (FL, USA) and Guwahati (Assam, India) but also for Santa Barbara (CA,
USA)! This fact is well-known to HVAC engineers but not necessarily to
control engineers working on VES. The only exception is in cold and dry climates in which economizers are used to bring in large amount of outdoor air to take advantage of free
cooling, and the exhaust ducts are
appropriately designed\footnote{Most buildings are not: large amount of outdoor air will create excessive pressurization  and doors will refuse to close, sometimes leading to alarms.}.

What is the implication of the discussion above for flexibility characterization? In most applications of optimization, only the optimal solution $x^*$ matters, the optimal value $f^*:=f(x^*)$ - where $f$ is the objective - is irrelevant. This is the case for energy efficient control of HVAC systems; as long as the optimal control computed reduces energy use from nominal conditions without violating constraints, it does not matter if the predicted optimal power is accurate or not. However, in flexibility characterization the objective $f(\cdot)$ is strongly dependent on demand $p$, such as the error $\|\rgrid - p\|$ in~\eqref{eq:planning}. In that case the optimal value $f^*= \|\rgrid - p^*\|$ is equally important. If a model of electrical demand $M_p:u\to p$ that ignores humidity is used the solution $p^*$ may be highly inaccurate. Such simplifications are sometimes made in designing control algorithms for providing VES; e.g.,~\cite{LiuOptimalIoT:2021,vrettos2016scheduling}. While use of feedback will help reduce the impact of modeling error, it is not clear if feedback alone can correct the large error in the model's prediction due to ignoring humidity. The argument for not ignoring humidity is stronger for methods for flexibility characterization, since robustness due to feedback is lacking in open loop flexibility characterization. 

\emph{The upshot is that if source of flexibile demand is that used for cooling and dehumidification, the model used for computing demand flexibility needs to incorporate
  humidity, in general.} An exception to this rule is when there is adequate time scale separation, so that the electrical demand that is manipulated for VES does not affect the cooling and dehumidification process, such as high frequency change in demand from fan power variation~\cite{linbarmeymid:2015} and chiller power variation~\cite{SuDemonstration-I-STBE:2015} to provide frequency regulation.


\ifx 0
  \subsubsection{When humidity can be ignored in modeling and control?} \rd{remove?}
  The answer differs based on the type of system. Let us first consider small packaged air conditioning units used in residential or small commercial buildings. The thermostat control logic in these units closes the loop based on temperature measurements alone. The reason is that temperature sensors are cheaper and more accurate than humidity sensors~\cite{}. As long as a residential air conditioner's cooling capacity is \emph{not too large} compared to the cooling load, the compressor runs sufficiently long to meet temperature setpoint. This has the effect of reducing the humidity of the indoor air since it is recirculated over the cooling coil. If the unit is so powerful that the temperature drops quickly upon turning on which turns the unit off, moisture removal is inadequate. This is a common source of poor humidity control, since many people are under the mistaken belief that a higher capacity unit is always better than a lower capacity unit~\cite{}.  

It can be argued that in addition to the requirement that $\theta \in [\theta_{\min}, \theta_{\max}]$, we also require that the average run time of the compressor is maintained by a VES control algorithm, humidity will be maintained to the same degree that the baseline thermostat logic.  This additional constraint can be easily incorporated into the framework of the definition~\eqref{}. Thus, humidity considerations can be ignored in modelign and control design for VES with small air conditioning units if a constraint is imposed on the minimal run time of the compressor in addition to temperature maintenance. If the unit is appropriately sized, maintaining indoor temperature and minimal run time will be not conflict. However, experimental work is needed to fully addressed the question: perhaps a minimal run time constraint is unnecessarily restrictive.  

The second part of the answer is for large commercial HVAC systems. If the time scale of the demand variation is much faster than the time scale of dehumidification in the cooling coil, we believe there will be no humidity violation. This is the situation when only fast VES services, such as frequency regulation, is provided by varying the fan speed~\cite{} or chiller XXX~\cite{}. Slower time scale services - such as load following or peak demand reduction - may not enjoy this advantage. In those case the effect of control commands on humidity need to be taken into account.

\subsubsection{Does humidity change the conclusions in Section~\ref{}?}
No???, the conclusions remain the same since the dynamics $p \to q$, from electrical power consumption to \emph{indoor humidity} is also low pass, so the same argument provided for temperature still applies. The low pass characteristics come from two sources. One, the dynamics $M_1$ from electrical demand to dehumidification of conditioned air at the AHU - which involves chillers and cooling
coils - is less than 2-3 minutes; see~\cite{}. Two, the dynamics $M_2$ from conditioned air flow rate to a zone to its space humidity ratio is low pass. Consider the following approximate model of $M_1$, which can be derived by a simple mass balance and assuming instant mixing of conditioned air delivered to the room with existing air:
\begin{align}
  \label{eq:8}
  \frac{dW}{dt} =  - \frac{1}{\rho_{air}V_{z}}\left(\dot{g}_{H_20}(t) + \mca\Wca(t) - \mca W(t)\right)
\end{align}
where $V_z$ is the zone volume, and $\dot{g}_{H_20}$ is the moisture generation rate insize the zone. The time constant of this first order system is $\tau = \rho_{air}V_{z}/\mca$. For an office room $12' \times 12' \times 10'$, with 1500 cfm air flow rate, the time constant turns out to be \rd{XXX} hour. Thus, the overall time constant from $p$ to $W$ is around $XXXXXXXX$.

There is a transport delay from chiller plant power consumption to indoor temperature and humidity in a building due to the flow of chilled water from the chiller at the plant to the cooling coil at the building's AHU. This delay will cause a challenge in tracking a power reference signal by manipulating chilled water flow rate at the AHU either directly, or indirectly by global thermostat set point change~\cite{}, or by manipulating airflow rate in the AHU that changes chilled water return temperature~\cite{}. Any VES control algorithm needs to take this delay into account~\cite{linbarmat:2017}. The often mentioned ``rebound effect'' seen in demand response attempts using thermostat set point changes most likely occurs because of a failure to do so.

\fi

\subsection{Flexibility characterization of large commercial HVAC systems is an opportunity}\label{sec:largeHVAC}
Hydronic HVAC systems - that use chilled water for cooling and dehumidification, or hot water for heating  - are commonly used in large buildings. Electrical power is consumed at multiple equipment and each has a different type of flexibility, which requires a careful bookkeeping. 

\begin{figure}[ht]
	\includegraphics[width=1.0\linewidth]{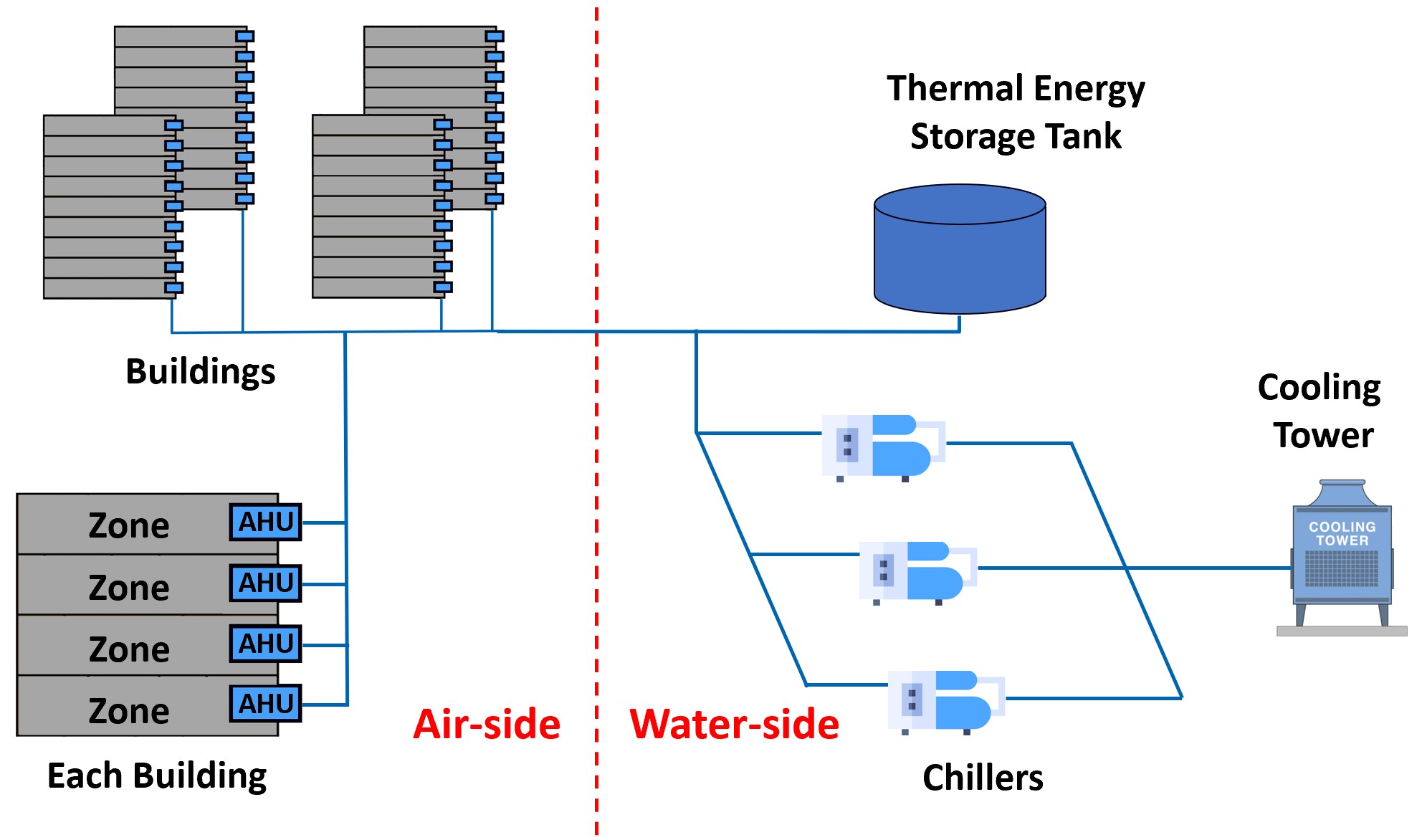}  
	\caption{A hydronic HVAC system used in large buildings. The TES shown is not common.}
	\label{fig:HVAC-plant}
      \end{figure}

Figure~\ref{fig:HVAC-plant} shows a commercial hydronic HVAC system that is common in the USA and in South East Asia. The chillers cool the warmer chilled water returning from the buildings, and the heat so extracted is rejected to the environment at the cooling towers. This ``water side'' equipment - chillers, pumps, cooling towers, are located in a ``chiller plant'' that is separate from the buildings they serve. The colder chilled water supplied to the buildings are used in one or more air handling units (AHUs) to cool and dehumidify air before supplied to the building's zones. An AHU is shown in Figure~\ref{fig:HVAC-single-zone}. Some chiller plants also have a thermal energy storage (TES) system, which stores either chilled water or ice. In colder climates, the water side equipment has boilers instead of chillers and the TES stores hot water.

\def\ppump{\ensuremath{p^{\text{pump}}}}
The total (electrical) power demand of such a system with chilled water consists of three main components: (i) power consumed in cooling and dehumidifying the mixed
air stream, $\pcc$, (ii) power consumed in reheating,
$\prh$, and (iii) power needed to circulate the air
$\pfan$, consumed by fan motors, and that required to circulate
water, $\ppump$, consumed by pump motors. It should be emphasized
that AHUs where cooling and dehumidification occurs, consume ``cooling'' from the chilled water loop, but the electricity to provide that cooling is consumed in
the chillers and the cooling towers, and there is a transport delay between the two. 

The last two components - reheating and fan/pump power - are simple to model; see~\cite{RamanAnalysisTSG:2020,Vrettos_etal_part1:2018,TianRealtimeTSG:2022}. But modeling electric demand for cooling and dehumidification is far more challenging, which involve both water-side and air-side processes. Models that relate control commands at the chiller plant (such as chilled water supply temperature set points and cooling tower fan speed set points at the supervisory level, or chiller on/off commands and chiller guide vane position commands at the lower level) to both the electrical demand at the chiller plant and the chilled water inlet conditions at the AHUs, are quite complex: first principles based models typically involve partial differential equations. Apart from strong nonlinearities, there is a transport delay between the electrical power consumption at the chilled water plant and the cooling consumed at the building. At the air side, control command and disturbances at the cooling coil (such as mixed air flow rate, temperature and humidity) to outputs (such as conditioned air humidity) and indoor climate conditions that determine QoS, are also quite challenging to model. It is not possible in this paper to review the large body of existing work on modeling these systems, or even to discuss the reasons for complexity. We refer the interested reader to~\cite{RamanMPC_AE:2020,GuoOptimalJESBC:2024,TianRealtimeTSG:2022} and references therein for an incomplete list of relevant work.

The only successful use of large commercial HVAC equipment for VES service so far has relied on time scale separation. In fact these successes - whether by using fans~\cite{linbarmeymid:2015,Vrettos_etal_part1:2018} or chillers~\cite{SuDemonstration-I-STBE:2015} - are all on frequency regulation, a fast service that has a time scale of less than a minute. By using time scale separation these works have avoided the need for complex equipment models.

However, we believe commercial HVAC systems have the potential to vary its demand in a longer time scale beyond frequency regulation since large buildings have large thermal inertia and the mass the chilled or hot water in the network of supply and return pipes provide additional inertia. This is especially true for HVAC systems with thermal energy storage. How long a time-scale is feasible is not yet known. Very little work has been in quantifying the flexibility potential of such HVAC systems, which is a missed opportunity. Complex dynamics of equipment and transport delay (of chilled water from the plant to the buildings) is a challenge in determining flexibility.




\section{Flexibility of a load collection (consumer focus)}\label{sec:many}
Since each load's demand flexibility is typically much smaller than what the grid needs as storage service, many loads will have to be coordinated to provide VES in real time. Perhaps as a result, both coordination algorithm design and demand flexibility characterization of load ensembles - rather than a single load - have been studied extensively.

The demand flexibility of a collection of $n$ loads can be defined as
\begin{align}
  \label{eq:6}
  \Omega = \{ \sum_i p_i | p_i \in \Omega_i, i = 1,\dots, n\}
\end{align}
where $\Omega_i$ is the flexibility set of load $i$ defined in~\eqref{eq:flexCapacity-general}.

Methods to characterize flexibility of the load collection should be (i) independent of the coordination algorithm. If they are dependent on a specific coordination algorithm, an improved coordination algorithm may lead to an increase in the capacity of the loads! 
These methods should also be (ii) computationally tractable.

Flexibility characterization of a collection of on/off thermostatically controlled loads (TCLs) has been the most extensively studied sub-topic under demand flexibility characterization. The on/off nature of individual loads make the problem challenging, which has led to many types of approximations; see~\cite{hao_aggregate:2015,ZhaoGeometricTPS:2017,BarotConciseIJEPES:2017,ZirasNewPSCC:2018,ChengHierarchicalEPES:2019,WangFlexibilityTSG:2020,coffmanAggregateTPS:2022} and references therein. Only some of these works provide methods satisfy the requirements (i) and (ii) specified in the previous paragraph, such as~\cite{hao_aggregate:2015,coffmanAggregateTPS:2022}. 

\begin{figure}[ht]
  \centering
  		\subfigure[junk][A load's flexibility (top) and how many loads can be coordinated to increase the collection's flexibility time scale (bottom). \label{fig:TCLcollection-increase-timescale-inexact}]
		{
			\includegraphics[width=0.6\columnwidth]{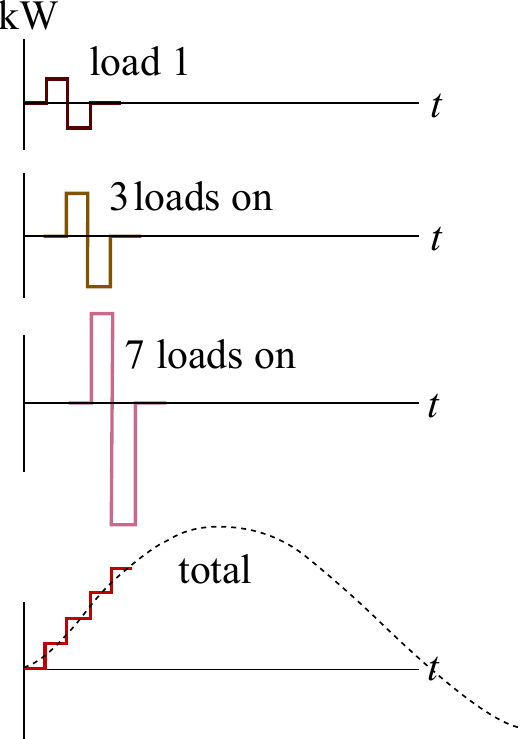}
		}
  		\subfigure[junk][A numerical example of increasing the flexibility time-scale through coordination. \label{fig:TCLcollection-increase-timescale-exact}]
		{
			\includegraphics[width=0.8\columnwidth]{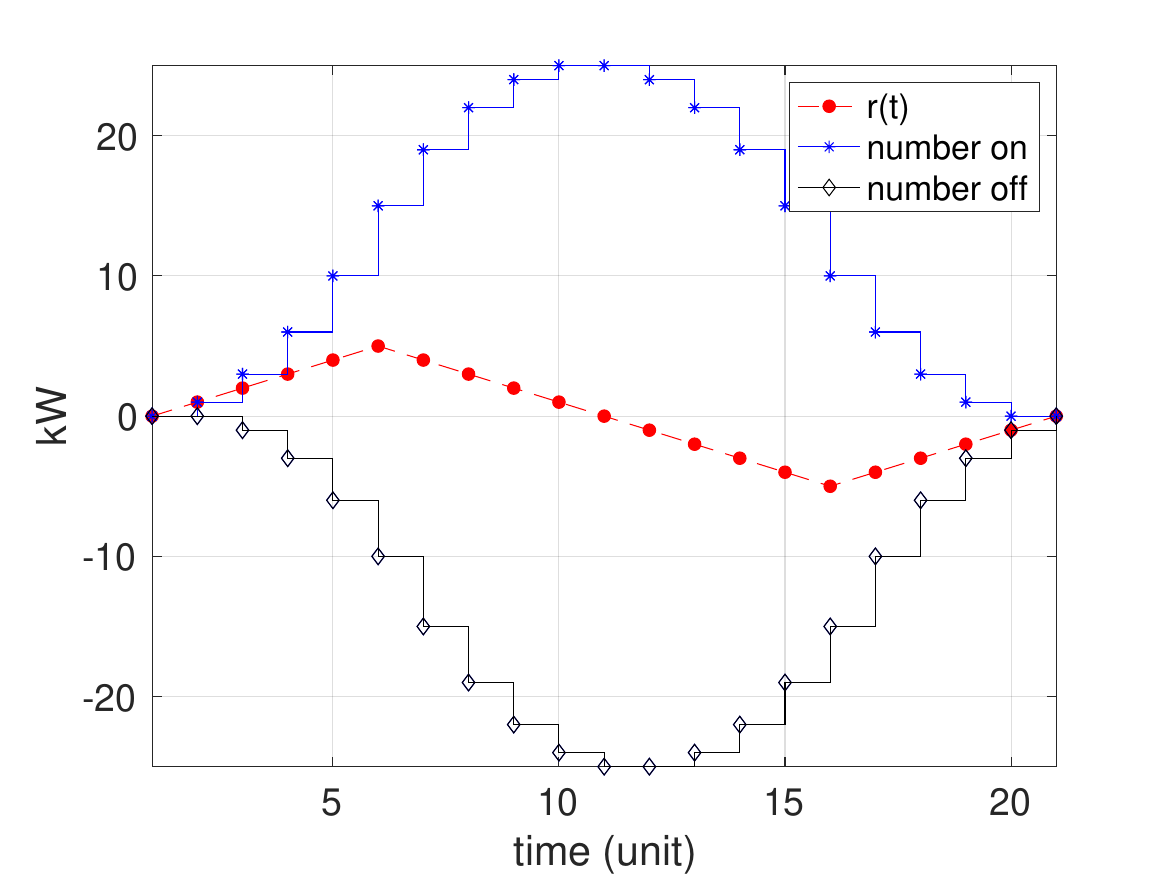}
		}
  \caption{Increase in timescale of virtual energy storage by coordinating a collection of loads.}
  \label{fig:TCLcollection-increase-timescale}
\end{figure}

It is obvious that if each load can deviate their demand by 1 kW, with $n$ loads the maximum possible deviation is $n$ kWs. What is not obvious - or at least has not been discussed much - is the qualitative difference in the \emph{time scale} of demand flexibility between an individual's and a collection's. To illustrate this difference, we consider a homogeneous collection in which each load's demand  can be varied in the discrete set $\{p^b,p^b+u,p^b-u\}$ where $p^b$ -in kW - is its (constant) baseline demand. Suppose each has a \emph{time-flexibility} of 1 time unit. Meaning, if the demand is increased from $p^b$ to $p^b+u$ kW, it must be held there for 1 time unit, and then the demand must be decreased to $p^b-u$ kW for the following time unit, to meet QoS constraints. This is shown in Figure~\ref{fig:TCLcollection-increase-timescale-inexact}(top). After the demand deviation becomes 0 at the end of two minutes, the cycle can be started again,
perhaps in opposite direction - first decrease in demand followed by
increase and then brigning back to the baseline. 
Imagine that the storage time scale required by the grid operator is
much larger than 1 time unit. An example of such a demand deviation reference is shown in Figure~\ref{fig:TCLcollection-increase-timescale-inexact} (bottom). An individual load cannot follow such a reference  no matter how large a flexibility amplitude $u$ is. However, it is
  possible to follow this reference closely with a large collection of loads even
  with a small individual flexibility.

  A specific example is shown in Figure~\ref{fig:TCLcollection-increase-timescale-exact}. Here the reference for demand deviation is a triangular waveform with period = 21 time units and peak amplitude 5 kW, while each load has a time flexibility of 1 time unit and demand flexibility of 1 kW. As the figure shows, this reference can tracked exactly with at most 24 loads. 

It is not too difficult to see that for a given size of the
collection, $n$, one can increase the peak amplitude of the storage
service but with a shorter time scale, and vice versa. In the
limit, the maximum amplitude will be $un$ at which the possible time
scale of the collection's VES service will be exactly the same as that of the
individual. At the other extreme, the maximum amplitude is only $u$ but the time-scale is much longer. 

The example above of the tradeoff between time scale and amplitude of
demand flexibility of a collection leads to a number of research
questions. For instance ``how to compute the 
  amplitude-vs-timescale curve for a given collection of loads''? And, ``how many air conditioners are needed to track a class of VES reference signals''? For a given collection of TCLs, and a given reference, the best the loads can do - meaning the reference closest to the given reference they can track without any load violating its own constraints - can be computed with the methods in~\cite{coffmanAggregateTPS:2022,CoffmanUnifiedAutomatica:2023}. Short of repeatedly trying these methods for an increasing number of loads, there is no technique available as yet that can provide the minimum number of loads needed to track a given reference. Another related question is, ``is there a specific timescale that is especially suitable for a particular type of load? Meaning, are water heaters better in renewable generation following than a collection of chiller plants?'' Ideas such as frequency-domain analysis of linearized models of load collections, as done in~\cite{meybarbusyueehr:TAC:2015}, can serve as a starting point in addressing the last question.

Whether for one load or for a collection, complexity of models is a hurdle in computing flexibility sets. So far, this hurdle is overcome by using low order approximations. But for some loads even low order models are hard to come by, such as large commercial HVAC systems with chillers and TES systems (see Sec.~\ref{sec:largeHVAC}). Similarly, degradation rates of EV batteries and arrival rates of EV charging demand from consumers are hard to model~\cite{LiLearning-BasedTSG:2021}. In contrast, sophisticated simulators of such systems exist, e.g. Modelica-based simulation models of complex HVAC systems~\cite{wetter2014modelica}. Methods of computing flexibility sets directly from data by using learning based techniques can be more valuable than methods that require low order models. This approach has been explored in~\cite{LiLearning-BasedTSG:2021,CoffmanModelfreeACC:2021}, but much more work is needed. 

\section{Grid operators' requirements on flexibility}\label{sec:industry}
In this section we discuss three issues related to demand flexibility that are of particular interest to balancing authorities and load aggregators.

\subsection{Battery equivalent characterization of flexible demand}
An electrochemical battery is typically characterized by its energy
capacity (Wh), and maximum charging and discharging rates (W). To aid
decision making by those in charge of choosing between real and
virtual batteries, \emph{these numbers need to be provided for VES as well}.

Often the maximum charging rate, the power capacity $p_{cap}$, is simply assumed to be the difference between the rated power of the equipment and the baseline power. However, such a demand deviation may not be maintainable at a specific time $t_0$ depending on what happened before $t_0$, as doing so may violate some QoS constraint such as temperature in an HVAC system. Similarly, a simple definition of energy capacity might be $p_{cap}
\times \tau$ where $\tau$ is the maximum duration that the power
deviation can be maintained at the power capacity $p_{cap}$. This definition too is problematic for the same reason, since it may not be possible to keep charging for $\tau$ duration without violating some QoS constraint depending on the initial QoS, $q(0)$. 

We now propose a formal definition. 
Recall that the power consumption of a virtual battery is $\tilde{p}(t):=p(t) -p^b(t)$, and is considered charging when $\tilde{p}(t)$ is positive and discharging when negative.  
The corresponding energy stored in a virtual battery at $t$, denoted by $\tilde{e}_{\tilde{p}}(t)$ is now defined as $\tilde{e}_{p}(t) := \int_{-\infty}^{t} \tilde{p}(v)dv$, which, unlike real batteries, can be negative as well. The maximum charging and discharging rates of the virtual battery are now defined as
\begin{align}
  p_c^w & = \max_{p \text{ feasible}}\{\max_t \tilde{p}(t)  \}, & p_{dc}^w & = \max_{p \text{ feasible}}\{\max_t (- \tilde{p}(t)) \} 
\end{align}
As before, a demand trajectory $p$ is called feasible if QoS is maintained by it, i.e., $q(t) \in Q(t), \forall t$ under $p$. The superscript $w$ denotes disturbance; capacities of a virtual battery depends on the associated disturbance trajectory. For an HVAC-based virtual battery, weather is the most significant part of disturbance. For instance, on a 
hot day an air conditioner has to operate near its rated power most of
the time under baseline conditions, and so its charging rate
$p_{c}^w$ will be small. The \emph{charging and discharging energy capacity} of the virtual battery are now defined as
\begin{align}
  e_c^w & = \max_{p \text{ feasible}}\{\max_t \tilde{e}(t)  \}, & e_{dc}^w & = \max_{p \text{ feasible}}\{\max_t (- \tilde{e}(t)) \} 
\end{align}


\emph{Computing these power and energy capacities for any type of HVAC loads is an open problem.} The paper~\cite{CoffmanCharacterizingTPS:2020} provides a method for a related capacity definition, but does
not take into account the effect of weather. The paper~\cite{HurtadoQuantifyingAE:2017} computes power and energy capacities with a different definition that is inspired by terminology from generators such as ramp rate and ramp duration. 

\subsection{Cost of demand flexibility service}
In our discussions with utilities and BAs that are looking for energy storage technologies, lack of knowledge of the cost of VES is frequently cited as a big barrier in its adoption. Batteries, though expensive, have a clear advantage: their costs can be estimated far more easily. Very little work has been done on estimating the cost of VES. A notable exception is~\cite{cammardella2018energy}, which computed the net present value of VES service of a collection of water heaters. The study concluded that in some cases, VES can be much cheaper than an electrochemical battery. But much more work is needed, for instance, to verify that the cost of large scale deployment assumed in~\cite{cammardella2018energy} is not unduly optimistic. 

Although a  large number of demand-flexibility demonstrations have been conducted with HVAC systems - see \cite{WaseemTechnologiesJMPSCE:2021} for an incomplete list - these demonstations use retrofitted appliances and thus do not inform us on the cost of large scale deployment in which loads will come with ``grid-friendly'' technologies installed at the factory. 
 
\subsection{Flexibility characterization for long term grid planning}
The flexibility discussed so far, $\Omega_{t_0}^{t_f}$ is for short term flexibility, with $t_f - t_0$ being of the order of a few minutes to a day, but not months or years. Consider an ancillary service market in which flexible loads or load aggregators participate. Uncertainty - due to weather etc. - will create an uncertainty in the flexibility capacity estimate, which will affect how flexible loads bid and how the market operators plan storage requirements. As long as the forecast horizon is short, say, less than a day, this uncertainty will be low since weather forecast is quite accurate in this short horizon. 

But for long term planning over months and years, the uncertainty in
forecasts is simply too great to be useful for both grid planners and loads providing VES. \emph{There is no consensus yet on how to define long term VES flexibility, much less compute it. } A statistical viewpoint, using spectral density, is proposed in~\cite{CoffmanCharacterizingTPS:2020}. The argument in~\cite{CoffmanCharacterizingTPS:2020} was that long term statistics of weather and grid conditions can be forecasted more accurately than time domain signals. But the bounds for load collections obtained were highly conservative. An answer to the long term flexibility question will help determine what fraction of a grid's requirements next year can be met by, say, 100,000 air conditioners.

\section{Conclusion}\label{sec:conclusion}
Some of the gaps in the literature on flexibility characterization might be due to the overemphasis on coordination algorithm design for load collections, which sometimes conflates the effect of the algorithm from the effect of constraints due to loads' properties and consumers' preferences. Characterizing load flexibility capacity that is independent of the coordination algorithm is paramount to avoid this confusion. Perhaps designing of centralized controllers to exploit flexible demand in the service of the power grid - even though they may not be suitable for deployment - should be explored to provide estimates of load capacities  that can serve as upper bound on what is achievable with distributed coordination. 

Due to lack of space, we have avoided delving into the topic of uncertainty. Methods are needed to assess the impact of uncertainty on flexibility of demand, such as weather on HVAC demand and consumer behavior on EV demand.

\bibliographystyle{IEEEtran}
\input{v1.bbl}

\end{document}